\documentclass[11pt,a4paper]{article}
\usepackage{ifthen,latexsym,amssymb,amsmath,amsthm,bbm,fixmath}
\usepackage[shortlabels]{enumitem}

\newcommand{\C}[1]{{\protect\mathcal{#1}}}
\newcommand{\I}[1]{{\mathbbm #1}}
\newcommand{\V}[1]{\mathbold{#1}}

\newcommand{\e}{\varepsilon}
\newcommand{\floor}[1]{\lfloor #1\rfloor}
\renewcommand{\mid}{:}
\renewcommand{\ge}{\geqslant}
\renewcommand{\le}{\leqslant}
\newcommand{\Matrix}[2]{\left(\begin{array}{#1} #2\end{array}\right)}

\newcommand{\beq}[1]{\begin{equation}\label{#1}}
\newcommand{\eeq}{\end{equation}}

\newtheorem{theorem}{Theorem}
\newtheorem{lemma}[theorem]{Lemma}
\newtheorem{question}[theorem]{Question}
\newtheorem{proposition}[theorem]{Proposition}
\newtheorem{corollary}[theorem]{Corollary}

\newcommand{\bpf}[1][Proof.]{\begin{proof}[#1] }
\newcommand{\epf}{\end{proof}}

\newtheorem{claim}{Claim}[section]
\newcommand{\bcpf}{\bpf[Proof of Claim.]}
\newcommand{\ecpf}{\epf}

\usepackage{authblk}

\theoremstyle{remark}
\newtheorem{remark}[theorem]{Remark}

\newcommand{\SosFamily}{\mathcal{W}_S}

\newcounter{ct1}
\newcounter{ct2}

\newcommand\blfootnote[1]{
  \begingroup
  \renewcommand\thefootnote{}\footnote{#1}
  \addtocounter{footnote}{-1}
  \endgroup
}

\begin{document}

\newcommand{\pr}{\mathbb{P}}

\newcommand{\codeg}{\mathrm{codeg}}
\newcommand{\dd}{\,\mathrm{d}}
\newcommand{\supp}{\mathrm{supp}}
\newcommand{\h}[2]{t(#1,#2)}

\newcommand{\Pxy}{Z}
\newcommand{\cutnorm}[2]{\|#1-#2\|_\Box}

\newcommand{\falling}[2]{\left(#1\right)_{#2}}
\newcommand{\fallingtwo}[2]{#1_{#2}}

\newcommand{\Lo}[1]{\cite[#1]{Lovasz:lngl}}

\renewcommand{\chi}{\I 1}

\setcounter{ct1}{1}
\setcounter{ct2}{2}
\author[1]{Oliver Cooley$^\fnsymbol{ct1}$}
\author[2]{Mihyun Kang$^\fnsymbol{ct1}$}
\author[3]{Oleg Pikhurko$^\fnsymbol{ct2}$}
\affil[1]{Institute of Science and Technology Austria (ISTA)\\ Am Campus 1, 3400 Klosterneuburg, Austria, \texttt{oliver.cooley@ist.ac.at}}
\affil[2]{Institute of Discrete Mathematics\\ Graz University of Technology\\ Steyrergasse 30, 8010 Graz, Austria, \texttt{kang@math.tugraz.at}}
\affil[3]{Mathematics Institute and DIMAP\\
	University of Warwick\\
	Coventry CV4 7AL, UK, \texttt{O.Pikhurko@warwick.ac.uk}}

\title{On a question of Vera T.\ S\'os\\
about size forcing of graphons}

\maketitle

\begin{abstract}
	The \emph{$k$-sample} $\I G(k,W)$ from a graphon $W:[0,1]^2\to [0,1]$ is the random graph on $\{1,\dots,k\}$, where we sample $x_1,\dots,x_k\in [0,1]$ uniformly at random and make each pair $\{i,j\}\subseteq \{1,\dots,k\}$ an edge with probability $W(x_i,x_j)$, with all these choices being mutually independent. Let the random variable $X_k(W)$ be the number of edges in~$\I G(k,W)$.

	Vera T.\ S\'os asked in 2012 whether two graphons $U,W$ are necessarily weakly isomorphic if the random variables $X_k(U)$ and $X_k(W)$ have the same distribution for every integer $k\ge 2$.
	 This question when one of the graphons $W$ is a constant function was answered positively by Endre Cs\'oka and independently by Jacob Fox, Tomasz {\L}uczak and Vera T.\ S\'os.
	Here we investigate the question when $W$ is a 2-step graphon and prove that the answer is positive for a 3-dimensional family of such graphons. 

	We also present some related results.
\end{abstract}

\renewcommand{\thefootnote}{\fnsymbol{footnote}}

\section{Introduction}
\footnotetext[1]{Supported by Austrian Science Fund (FWF) Grant~I3747.}
\footnotetext[2]{Supported by ERC Advanced Grant 101020255 and Leverhulme Research Project Grant RPG-2018-424.}
\blfootnote{\emph{Key words and phrases:} Graphons, $k$-sample, forcing, containers.}
\blfootnote{\emph{Mathematics subject classification:} 05C99, 05C80}
\renewcommand{\thefootnote}{\arabic{footnote}}
\blfootnote{An extended abstract of this paper appeared in the Proceedings of the European Conference on Combinatorics, Graph Theory and Applications (EuroComb 2021), CRM Research Perspectives, Springer.}
\emph{Graphons} (that is, measurable symmetric functions $[0,1]^2\to [0,1]$) have recently found many important applications in other areas, such as the limit theory of dense graphs (Lov\'{a}sz et al~\cite{BCLSV06,LovaszSzegedy06,LovaszSzegedy07gafa}),
 large deviation principles for random graphs (Chatterjee and Varadhan~\cite{ChatterjeeVaradhan11}),
property testing in computer science (Lov\'{a}sz and Szegedy~\cite{LovaszSzegedy10ijm}), etc.  We refer the reader  to the monograph by Lov\'asz~\cite{Lovasz:lngl} for an introduction.

The \emph{$k$-sample} $\I G(k,W)$ from a graphon $W$ is the random graph on $[k]:=\{1,\dots,k\}$ obtained by  sampling $x_1,\dots,x_k\in [0,1]$ uniformly at random and making each pair $\{i,j\}\subseteq [k]$ an edge with probability $W(x_i,x_j)$, with all these choices being mutually independent. The \emph{(homomorphism) density} $\h{F}{W}$ of a graph $F$ on $[k]$ in $W$ is the probability that $E(F)\subseteq E(\I G(k,W))$, that is, every adjacent pair in $F$ is also adjacent in~$\I G(k,W)$.
Equivalently, we can define
\beq{eq:t=Int}
 \h{F}{W}:=\int_{[0,1]^k} \prod_{\{i,j\}\in E(F)} W(x_i,x_j)\dd x_1\dots\dd x_k.
 \eeq
 Let us call two graphons $U$ and $W$ \emph{weakly isomorphic} if the random graphs $\I G(k,U)$ and $\I G(k,W)$ have the same distribution for every~$k\in\I N$. This is equivalent to 
$\h{H}{U}=\h{H}{W}$ for every connected graph~$H$, because we can recover the distribution of $\I G(k,W)$ using the following inclusion-exclusion formula
 \beq{eq:5.20}
 \I P(\I G(k,W)=G)=\sum_{F\supseteq G\atop V(F)=[k]} (-1)^{|E(F)\setminus E(G)|}\,\h{F}{W},\quad\mbox{for every graph $G$ on $[k]$,}
 \eeq
 and replacing each $\h{F}{W}$ by the product of the densities of the components of~$F$ (see Lemma~\ref{lm:ProductDensity} later).

Borgs, Chayes and Lov{\'a}sz~\cite{BorgsChayesLovasz10}
showed that all graphons in the weak isomorphism class of $W$  can, roughly speaking, be obtained from $W$ by applying measure-preserving transformations of the variables. (See also Diaconis and Janson~\cite{DiaconisJanson08} who  derived this result from the Aldous--Hoover Theorem~\cite{Aldous81,Hoover79} by noting a connection to exchangeable arrays.) This gives an analogue of the classical moment problem, where each $\h{F}{W}$ can be thought of as the ``$F$-th moment'' of~$W$.

A \emph{graphon parameter} $f$ is a function that assigns to each graphon $W$ a real number or a real vector $f(W)$ such that $f(W)=f(U)$ whenever $U$ and $W$ are weakly isomorphic. We say that a family $(f_i)_{i\in I}$ of graphon parameters \emph{forces} a graphon $W$ if every graphon $U$ with $f_i(U)=f_i(W)$ for every $i\in I$ is weakly isomorphic to~$W$. For example, the famous result of Chung, Graham and Wilson~\cite{ChungGrahamWilson89} on $p$-quasirandom graphs can be stated in this language as follows.

\begin{theorem}[Chung, Graham and Wilson~\cite{ChungGrahamWilson89}]
\label{th:CGW} The constant-$p$ graphon is forced by $\h{K_2}{\cdot}$ and $\h{C_4}{\cdot}$, that is, by  the edge and 4-cycle densities.\qed\end{theorem}

Call a family $(f_i)_{i\in I}$ of graphon parameters \emph{forcing} if it forces every graphon~$W$.  For example, the densities $\h{F}{\cdot}$, where $F$ ranges over all
connected
 graphs, form a forcing family (by 
 \eqref{eq:5.20} and Lemma~\ref{lm:ProductDensity}). 

A lot of effort has gone into investigating whether a graphon $W$ is forced   by much less information than the densities of all graphs. Graphons that are forced by a finite set of graph densities are called \emph{finitely forcible} and their systematic study was initiated by Lov{\'a}sz, S{\'o}s and Szegedy~\cite{LovaszSos08,LovaszSzegedy11},
motivated by quasirandom graphs and extremal graph theory. As one would expect, finitely forcible graphons are ``rare'': they form a meagre subset of the space of all graphons (\cite[Theorem~7.12]{LovaszSzegedy11}).

The authors are not aware of any results where a substantially smaller set of parameters than the densities of all connected graphs is shown to be forcing. Vera T.\ S\'os~\cite{Sos12} posed some questions in this direction, and in particular considered the following problem.
For a graphon $W$ and an integer $k\in\I N$, let $X_k(W):=|E(\I G(k,W))|$ be the size of, i.e.\ number of edges in, the $k$-sample $\I G(k,W)$ from~$W$. We identify the random variable $X_k(W)$ with the vector of probabilities $\I P(X_k(W)=i)$ for ${0\le i\le {k\choose 2}}$, viewing it as a graphon parameter.
	Let $\SosFamily$ be the family of graphons $W$ that are forced by the sequence $(X_k(W))_{k\in\I N}$, i.e.\ by the distributions of sizes of samples from $W$.

\begin{question}[Size Forcing Question (S\'os~\cite{Sos12})]
	\label{q:Sos}
	Is every graphon in $\SosFamily$?
\end{question}

Noga Alon (unpublished, see~\cite{Csoka16}) and independently Jakub Sliacan~\cite{Sliacan15} proved that the constant~$\frac12$ graphon is in the family $\SosFamily$. Then Endre Cs\'oka~\cite{Csoka16} and independently Jacob Fox, Tomasz {\L}uczak and Vera T.\ S\'os~\cite{FoxPersonal}
proved that constant-$p$ graphon is in the family $\SosFamily$ for any $p\in(0,1)$. A natural next step would be to try to
determine whether $W\in \SosFamily$ when $W$ is a \emph{2-step graphon}, that is, we have a partition of $[0,1]$ into two measurable sets $A$ and $B$ such that $W$ is constant on each of the sets $A^2$, $B^2$ and $(A\times B)\cup (B\times A)$. 
By replacing $W$ by a weakly isomorphic graphon, we can assume that $A=[0,a)$ and $B=[a,1]$ are intervals. Thus we need four parameters to describe a 2-step graphon: the measure of $A$ as well as the three possible values of $W$. 

Unfortunately, we were not able to prove that $W\in \SosFamily$ for every 2-step graphon~$W$.
However, we could prove this for the following 3-dimensional set of graphons.

\begin{theorem}\label{th:Negated} Let $W$ be the 2-step graphon with parts $A:=[0,a)$ and $B:=[a,1]$ such that its values on $A^2$, $(A\times B)\cup (B\times A)$ and $B^2$ are respectively $0$, $p\in (0,1]$ and~$q\in (0,1]$.
	If  $(1-a)q\le (1-2a)p$, then $W \in \SosFamily$.
	\end{theorem}
Let us mention here that since $X_k(1-W)$ has the same distribution as $\binom{k}{2}-X_k(W)$ (by taking complements),
if $W$ is forced by some sub-family of $(X_k)_{k \in \I N}$, then so is $1-W$.
Thus if $W$ lies in $\SosFamily$, then so does $1-W$.

We can also answer Question~\ref{q:Sos} for some other families of 2-step graphons $W$. Here we present two further examples (Theorems~\ref{th:1param} and~\ref{th:01p}) where a \emph{finite} set of some natural real-valued parameters suffices.

The first is motivated by the result of Cs\'oka~\cite{Csoka16} who in fact proved that the constant-$p$ graphon is forced by $X_4$ alone.
The following theorem proves a similar claim to that of Cs\'oka~\cite{Csoka16} for the limit of balanced quasirandom bipartite graphs, namely that it is forced by the edge distribution of its $5$-sample.

\begin{theorem}\label{th:1param}
	Let $p\in [0,1]$ and let $W$ be the graphon which is $0$ on $[0,1/2)^2\cup [1/2,1]^2$ and $p$ everywhere else. 
	Then $W$ is forced by $X_5$ alone.
\end{theorem}

Let the \emph{independence ratio} $\alpha(W)$ of a graphon $W$ be the supremum of the measure of $A\subseteq [0,1]$ such that $W(x,y)=0$ for a.e.~$(x,y)\in A^2$. As was observed by Hladk\'y, Hu and Piguet~\cite[Lemma~2.4]{HladkyHuPiguet19}, the supremum is in fact a maximum (that is, it is attained by some~$A$). Also, the \emph{clique ratio} $\omega(W):=\alpha(1-W)$ is the maximum measure of $A\subseteq [0,1]$ with $W$ being 1 a.e.\ on~$A^2$.

\begin{theorem}\label{th:01p}
	Given $a,p\in [0,1]$, 
	set $A:=[0,a)$ and $B:=[a,1]$, and
	let $W$ be the graphon which is 0 on $A^2$, $1$ on $B^2$, and $p$ everywhere else.
	Then $W$ is forced by $(\alpha,\omega,X_4)$.\end{theorem}

By using a basic version of the container method, we show that the value of $\alpha$ (and thus of~$\omega$) is determined by any infinite subsequence of $(X_k)_{k\in\I N}$. More precisely, the following holds.

\begin{theorem}\label{th:alpha}
	For every graphon $W$, it holds that 
	$$
	\alpha(W)=\lim_{k\to\infty} \big(\I P(X_k(W)=0)\big)^{1/k}.
	$$
\end{theorem}

Hladk\'y and Rocha~\cite{HladkyRocha20} defined and studied graphon versions of various graph parameters, including the independence ratio $\alpha(W)$. In particular, they investigated how these parameters 
can be related to
graph densities. Our Theorem~\ref{th:alpha}, by relating $\alpha(W)$ to graph densities, fills one missing entry in~\cite[Table~1]{HladkyRocha20}.

By combining Theorems~\ref{th:01p} and~\ref{th:alpha}, we directly obtain the following result.

\begin{corollary} Let $W$ be a graphon as in Theorem~\ref{th:01p} (that is, $W$ is
	$0$ on $[0,a)^2$, $1$ on $[a,1]^2$, and $p$ everywhere else). Then $W \in \SosFamily$.\qed
	\end{corollary}

Call a family $\C F$ of graphs \emph{forcing} if the corresponding family of parameters $(\h{F}{\cdot})_{F\in\C F}$ is forcing. 
S\'os~\cite{Sos12} also asked
if one can find substantially smaller forcing families than taking all connected graphs. We show that two natural examples,
namely the family of all cycles and the family of all complete bipartite graphs, do not suffice.

\begin{proposition}\label{pr:cycles} The family of all connected graphs with at most one cycle is not forcing. In particular, the family of all cycles is not forcing.\end{proposition}

\begin{proposition}\label{pr:Kkl} For every integer $d$, the family of all graphs of diameter at most $d$ is not forcing. In particular, the family of all complete bipartite graphs is not forcing.\end{proposition}

Also, let us mention here the somewhat related result of Shapira and Tyomkin~\cite{ShapiraTyomkin21} that the constant-$p$ graphon is not forced by $(\h{K_k}{\cdot})_{k\in\I N}$, that is, by clique densities.

\subsubsection*{Paper overview}

The paper is arranged as follows. In Section~\ref{Notation} we recall various standard notation that we will use, in particular notation related to graphons.
In Section~\ref{aux} we collect some easy preliminary results which we will apply later.
Theorem~\ref{th:1param} is then proved in Section~\ref{sec:1param}.

In Section~\ref{containers}, we first present an auxiliary graph result (Theorem~\ref{th:ManyEk}) which relates small independent sets to large almost independent sets, and prove this result using the container method. We subsequently use this result to prove Theorem~\ref{th:alpha}.

Sections~\ref{sec:Negated} and~\ref{01p} contain the proofs of Theorems~\ref{th:Negated} and~\ref{th:01p}, respectively.
Finally, Propositions~\ref{pr:cycles} and~\ref{pr:Kkl} are proved in Section~\ref{OtherQns}.
Section~\ref{sec:concluding} contains concluding remarks.

\section{Notation}\label{Notation}

Here we present some notation that is used in this paper.

Let $k$ be a non-negative integer. We denote the \emph{falling factorial} of a real number
$r$ by
$$
\falling{r}{k}:=r(r-1)\ldots (r-k+1).
$$
For a set $X$, let 
$${X\choose k}:=\{Y\subseteq X\mid |Y|=k\}$$
consist of all subsets of $X$ of size~$k$. 
The \emph{characteristic function} $\chi_X$ of $X$ assumes value 1 on $X$ and 0 everywhere else. We may abbreviate an unordered pair $\{x,y\}$ to~$xy$.

We will use the following notation related to graphs.
Let $G=(V,E)$ be a graph. Its \emph{complement} is~$\overline{G}:=\big(V,{V\choose 2}\setminus E\big)$. For $A\subseteq V(G)$, its \emph{neighbourhood}
$$
N(A):=\{x\mid \exists\, y\in A\mbox{ with } xy\in E\}
$$ 
 consists of vertices that send at least one edge to~$A$. For graphs $H_1$ and $H_2$, their \emph{disjoint union} $H_1\sqcup H_2$ is obtained by taking the union of vertex-disjoint copies of these graphs (with no edges across).

We will also be using the following special graphs. 
The \emph{$k$-clique} $K_k$ is the graph on $[k]$ in which every two vertices are adjacent. 
The \emph{$k$-path} $P_k$ is the path on $[k]$ that visits vertices $1,\dots,k$ in this order. 
The \emph{$k$-cycle} $C_k$ is the cycle on $[k]$ that visits vertices $1,\dots,k$ in this cyclic order.
Also, $K_{k,\ell}$ denotes the complete bipartite graph with parts of sizes $k$ and~$\ell$.

Let $\C G_{k,m}$ consist of isomorphism classes of all graphs with  at most $k$ vertices and exactly $m$ edges that do not contain any isolated vertices. For example, $\C G_{5,3}=\{K_3,P_4,P_3\sqcup K_2,K_{1,3}\}$.

Let us also collect some definitions related to graphons. 

The unit interval $[0,1]$ is by default equipped with the Lebesgue measure, denoted by~$\lambda$. When making any statements about subsets of $[0,1]$, we usually mean that they hold up to a set of measure~$0$. 
We will often use (Fubini-)Tonelli's theorem (see e.g.\ \cite[Theorem 14.2]{Dibenedetto16ra}) whose main part states, informally speaking, that non-negative measurable functions can be integrated in any order of variables. In particular, when working with $\h{F}{W}$ as the value of the integral in~\eqref{eq:t=Int}, we can integrate the variables $x_1,\dots,x_k$ in any order.
We will therefore often change the order of integration without mentioning this theorem explicitly.

Let $W$ be a graphon and let $A\subseteq [0,1]$ be a measurable subset. The \emph{degree} (resp.\ \emph{$A$-degree}) of $x\in [0,1]$ is 
$$\deg^W(x):=\int_0^1 W(x,y)\dd y$$
 (resp.\
$\deg^W_A(x):=\int_A W(x,y)\dd y$). The degree is defined for a.e.\ $x\in [0,1]$ by a part of Tonelli's theorem.
We call $W$ \emph{$p$-regular} if $\deg^W(x)=p$ for a.e.~$x\in [0,1]$.
The \emph{codegree} (resp.\ \emph{$A$-codegree}) of $(x,y)\in [0,1]^2$ is 
$$
\codeg^W(x,y):=\int_0^1 W(x,z)W(z,y)\dd z
$$
 (resp.\ $\codeg_A^W(x,y):=\int_A W(x,z)W(z,y)\dd z$). One can view $\codeg^W(x,y)$ as the density of 2-edge paths that connect $x$ and~$y$.

When discussing $\I G(k,W)$, it will often be convenient to view it as a graph whose vertex set consists of the sampled points $x_1,\dots,x_k\in [0,1]$ (which are  pairwise distinct with probability~1). Thus a phrase like ``$x_i$ is adjacent to $x_j$'' will be a shorthand for  ``$i$ is adjacent to $j$ in $\I G(k,W)$'', etc.

For a graph $G$ on $[k]$, its graphon $W_G$ is the graphon which is $1$ on $[\frac {i-1}k,\frac{i}k)\times [\frac{j-1}k,\frac jk)$ for each edge $ij\in E(G)$, and 0 everywhere else. In other words, we partition $[0,1)$ into $k$ intervals of length $1/k$ and let $W_G$ be the $k$-step $\{0,1\}$-valued graphon that naturally encodes the adjacency relation of~$G$.

\section{Some auxiliary results}\label{aux}

Let us present here  some known or easy auxiliary results that we need in this paper.

 \begin{lemma}\label{lm:ProductDensity}
 	For any graphon $W$ and for any graphs $H_1$ and $H_2$, we have
 	$$
 	\h{H_1 \sqcup H_2}{W}= \h{H_1}{W}\,\h{H_2}{W}.
 	$$
 \end{lemma}
 
 \bpf 
 Assume that $V(H_1\sqcup H_2)=[k]$ and let $A_1\cup A_2=[k]$ be the partition 
 into the vertex sets of $H_1$ and $H_2$. The lemma follows by observing that the subgraphs induced by $A_1$ and $A_2$ in $\I G(k,W)$ are independent of each other and, up to a relabelling of vertices, are distributed as $\I G(|A_1|,W)$ and  $\I G(|A_2|,W)$.
 \epf

 \begin{lemma}\label{lm:regular}
 	If $W$ is a $p$-regular graphon and $F'$ is obtained from a graph $F$ by attaching a pendant edge 
 	then $$\h{F'}{W}=p\,\h{F}{W}.$$
 \end{lemma}	
 
 \bpf We can assume that $V(F)=[k]$ and that the added edge is $\{k,k+1\}$. When computing $\h{F'}{W}$ as the integral over $(x_1,\dots,x_{k+1})\in [0,1]^{k+1}$ as in~\eqref{eq:t=Int}, we can 
 first integrate over~$x_{k+1}$. The only factor that depends on $x_{k+1}$ is $W(x_k,x_{k+1})$. Its integral is $p$ for a.e.\ $x_k$, so integrating out $x_{k+1}$ amounts to multiplying by~$p$ (and replacing $F'$ by $F$ in~\eqref{eq:t=Int}), proving the lemma.
 \epf

The following result implicitly appears in Cs\'oka~\cite{Csoka16}. For completeness, we present its proof.

\begin{lemma}\label{lm:Csoka} Let integers $k$ and $m$ satisfy $1\le m\le {k\choose 2}$. 
Then for every graphon $W$ we have
\beq{eq:Csoka}
\I E\left(\falling{X_k(W)}{m}\right)=\sum_{F\in \C G_{k,m}} c_{k,F}\,
 \h{F}{W},
\eeq
 where $c_{k,F}>0$ is $m!$ times the number of graphs on $[k]$ that, after discarding isolated vertices, are isomorphic to~$F$.
\end{lemma}

\bpf 
Let $\C X$ consist of all ordered $m$-tuples $(\{s_i,t_i\})_{i=1}^m$ of pairwise distinct pairs from~${[k]\choose 2}$. Thus, for example, its size $|\C X|$ is the falling factorial~$\fallingtwo{\big({k\choose 2}\big)}{m}$.   For $F\in\C G_{k,m}$, let $\C X_F$ consist of those sequences in $\C X$ that give a graph isomorphic to~$F$ after we discard all isolated vertices. Clearly, the sets $\C X_F$ when $F$ ranges over $\C G_{k,m}$ partition~$\C X$.

The left-hand side of~\eqref{eq:Csoka} is the expectation of the number of sequences in $\C X$ all of whose $m$ pairs are edges when we take the $k$-sample~$\I G(k,W)$. This expectation can be written as the sum over all  $(\{s_i,t_i\})_{i=1}^m\in \C X$ of the probability that each $\{s_i,t_i\}$ is an edge. The last probability is exactly $\h{F}{W}$ where $F$ is the unique graph of $\C G_{k,m}$ with $(\{s_i,t_i\})_{i=1}^m\in \C X_F$. The lemma follows by observing that each $F\in \C G_{k,m}$ appears exactly 
$|\C X_{F}|=c_{k,F}$ times this way.\epf

Here is a useful consequence of this lemma.

\begin{lemma}\label{lm:2Edge} Let $U$ and $W$ be graphons such that $X_k(U)$ and $X_k(W)$ have the same distributions for some $k\ge 3$. Then $\h{K_2}{U}=\h{K_2}{W}$, $\h{K_2\sqcup K_2}{U}=\h{K_2\sqcup K_2}{W}$ and $\h{P_3}{U}=\h{P_3}{W}$.\end{lemma}

\bpf By applying Lemma~\ref{lm:Csoka} with $m=1$ to $U$ and $W$, we get that 
$$
 {k\choose 2}\h{K_2}{U}=\I E\left(X_k(U)\right)=\I E\left(X_k(W)\right)={k\choose 2}\h{K_2}{W}.
$$
 Thus $U$ and $W$ have the same density of $K_2$ and, by Lemma~\ref{lm:ProductDensity}, of $K_2\sqcup K_2$. If $k\ge 4$, then Lemma~\ref{lm:Csoka} with $m=2$ gives that
 \begin{eqnarray*}
2! \, \frac{\falling{k}{3}}{2}\cdot \h{P_3}{U}&=&\I E\left(\falling{X_k(U)}{2}\right)-2!\,\frac{\falling{k}{4}}{8}\cdot \h{K_2\sqcup K_2}{U}\\
 &=&\I E\left(\falling{X_k(W)}{2}\right)-2! \,\frac{\falling{k}{4}}{8}\cdot \h{K_2\sqcup K_2}{W}
 \ =\ 2!\,\frac{\falling{k}{3}}{2}\cdot  \h{P_3}{W},
 \end{eqnarray*}
 finishing the proof (since $\falling{k}{3}\neq 0$). The same calculation applies for $k=3$ except that the $K_2\sqcup K_2$ term is absent.
\epf

We will also need the following bipartite analogue of Theorem~\ref{th:CGW}. While a
rather elementary proof by passing to finite graphs that converge to $U$ is possible (along the same lines as the original proof of Chung, Graham and Wilson~\cite{ChungGrahamWilson89}, see also e.g.\ \Lo{Theorem 11.62}), we present a proof that, while requiring some analytic background, deals directly with graphons.
	
\begin{lemma}\label{lm:BipQR} Let $A$ and $B$ be sets of measure $a$ and $b$ respectively that partition~$[0,1]$. (Thus $a+b=1$.) Let $p\in [0,1]$. Let $U$ be a graphon taking value 0 on $A^2\cup B^2$ such that 
	$\h{K_2}{U}=2abp$
	and $\h{C_4}{U}=2a^2b^2p^4$. Then $U(x,y)=p$ for a.e.\ $(x,y)\in (A\times B)\cup (B\times A)$.
\end{lemma}

\bpf
	Assume that $a,b\in(0,1)$ as otherwise there is nothing to do.
	
	Using that $U$ is 0 on $A^2\cup B^2$, we have 
	that $\h{C_4}{U}$ (that is, the density of the $4$-cycle in $U$) is $
	2\int_{A^2} \codeg^U_B(x,y)^2 \dd x\dd y$, where the factor 2 comes from having a partition into two equiprobable events, namely that $x_1,x_3\in A$ and that $x_2,x_4\in A$. On the other hand, by applying the Cauchy--Schwarz Inequality twice,
	we get that
	\begin{eqnarray*}
	2a^2b^2p^4 &=& 2\int_{A^2} \left(\codeg^U_B(x,y)\right)^2 \dd x\dd y\ \ge\ \frac2{a^2}\left(\int_{A^2} \codeg^U_B(x,y) \dd x\dd y
	\right)^2\\ &=& \frac2{a^2}\left(\int_B(\deg^U_A(z))^2\dd z\right)^2\
	\ge\ \frac2{a^2b^2}\left(\int_B\deg^U_A(z)\dd z\right)^4\  =\ 2a^2b^2p^4.
	\end{eqnarray*}
	Thus we have equality. This implies that $\codeg^U_B(x,y)= bp^2$ for a.e.\ $(x,y)\in A^2$ and that $\deg^U_A(x)=ap$ for a.e.~$x\in B$. The same argument applies when we count $C_4$ from the other side, giving that $\codeg^U_A(x,y)=ap^2$ for a.e.\ $(x,y)\in B^2$ and $\deg^U_B(x)=bp$ for a.e.~$x\in A$.

	View $U$ as the integral kernel operator 
	 defined by 
	 $$
	 (U\!f)(x):=\int_0^1 U(x,y)f(y)\dd y,\quad \mbox{for $f\in L^2([0,1])$ and $x\in [0,1]$.}
	 $$
	  The Cauchy-Schwarz (or H\"older's) Inequality gives that $U\!f\in L^2([0,1])$, so $U$ is an operator on~$L^2([0,1])$.
	  This operator is self-adjoint (since the function $U$ is symmetric) and compact (as an integral operator with its kernel $U$ being a bounded
	   and thus square-integrable function on~$[0,1]^2$, see e.g.\ \cite[Example 3 of Section 2.16]{GohbergGoldbergKaashoek03bclo}). The Spectral Decomposition Theorem (see e.g.\ 
	  \cite[Theorem 5.1 of Section 4.5]{GohbergGoldbergKaashoek03bclo}) 
	  gives that $L^2([0,1])$ has an orthonormal basis of eigenfunctions $(f_i)_{i\in\I N}$ with the corresponding eigenvalues $(\lambda_i)_{i\in \I N}$ such that $\lambda_i\to 0$ as $i\to\infty$. Then it follows that 
	  \beq{eq:SD}
	  U(x,y)=\sum_{i=1}^\infty \lambda_i f_i(x)f_i(y)\quad \mbox{for a.e.\ $(x,y)\in [0,1]^2$}.
	  \eeq

Consider the composition of $U$ with itself. This is again an integral kernel operator and we identify it with its kernel
	$$(U\circ U)(x,y):=\int U(x,z)U(z,y)\dd z=\codeg^U(x,y),\quad x,y\in [0,1].
	$$
	By above, $U\circ U$ is a.e.\ the two-step graphon of value 0 on $(A\times B)\cup (B\times A)$, value $bp^2$ on $A^2$ and value $ap^2$ on~$B^2$. (Recall that $U$ is 0 on $A^2\cup B^2$.) Thus, for every $f\in L^2([0,1])$, its image $(U\circ U)(f)$ is a function which is constant on $A$ and on~$B$. Thus, as an operator, $U\circ U$ has rank at most 2.  On the other hand, each of the characteristic functions $\chi_A$ and~$\chi_B$
	is an eigenvector of $U\circ U$ with eigenvalue~$abp^2$. We conclude that the operator $U\circ U$ has exactly one non-zero eigenvalue $abp^2$ of multiplicity~2.

  Clearly, the same functions $(f_i)_{i\in \I N}$ and the squares $(\lambda_i^2)_{i\in \I N}$ give a spectral decomposition of~$U\circ U$. 
Thus the rank of $U$ is also $2$ and its non-zero eigenvalues are in $\{-p \sqrt{ab},\,p \sqrt{ab}\,\}$. 
Using the established values for the (constant) degrees in $U$ across the partition $A\cup B$, we have that
$U\chi_A=ap \chi_B$ and $U\chi_B=bp\chi_A$.
Consider the functions $h_1:=\sqrt{b}\,\chi_A- \sqrt{a}\, \chi_B$ and $h_2:=\sqrt{b}\,\chi_A+ \sqrt{a}\,\chi_B$ that are orthogonal to each other and have $L^2$-norm~$\sqrt{2ab}$.
We have 
$$
 Uh_1=\sqrt{b}\, U\chi_A -\sqrt{a}\, U\chi_B  = ap\sqrt{b}\,\chi_B-bp\sqrt{a}\,\chi_A=-p\sqrt{ab}\, h_1,
 $$
 and similarly $Uh_2=p\sqrt{ab}\, h_2$. Thus, up to relabelling, we have $\lambda_1=-p\sqrt{ab}$ and $\lambda_2=p\sqrt{ab}$ (and $\lambda_i=0$ for all $i\ge 3$). Moreover, by the 1-dimensionality of the eigenspaces for $\lambda_1\not=\lambda_2$, it holds that $f_i=\pm h_i/\sqrt{2ab}$ for $i=1,2$. 
 By~\eqref{eq:SD}, we have that
	$$
	U(x,y)=\lambda_1f_1(x)f_1(y)+\lambda_2f_2(x)f_2(y) = \frac{-p\sqrt{ab}\,h_1(x)h_1(y)+p\sqrt{ab}\,h_2(x)h_2(y)}{{2ab}}.
	$$
	Thus $U$ is constant on $A\times B$ and its value there can be shown to be $p$
	by using the definition of $h_1$ and $h_2$ (or by our assumption that $\h{K_2}{U}=2abp$), giving the lemma.
\epf

\section{Proof of Theorem~\ref{th:1param}}\label{sec:1param}

Recall that $W$ is the 2-step graphon which is the limit of balanced bipartite $p$-quasirandom graphs.
Let $U$ be an arbitrary graphon such that the distribution of $X_5(U)$ is the same as the distribution of~$X_5(W)$. Let us denote this common distribution by~$X_5$. We will be iteratively proving a sequence of claims about $U$, until the derived information is enough to conclude that $U$ must be weakly isomorphic to~$W$. Assume that $p\not=0$, as the constant-0 graphon is clearly forced by $X_5$ being $0$ with probability~$1$.

\begin{claim}\label{prop:preg}
	The graphon $U$ is $(p/2)$-regular, that is, the degree function $\deg^U(x)=\int_0^1 U(x,y)\dd y$ is equal to $p/2$ for a.e.~$x\in [0,1]$.
\end{claim}

\bpf Consider the random variable $D:=\deg^U(x)$, where 
 $x\in [0,1]$ is uniformly random.  Lemma~\ref{lm:2Edge} shows that $\I E(D)=\h{K_2}{U}$ equals $\h{K_2}{W}=p/2$ and $\I E(D^2)=\h{P_3}{U}$ equals $\h{P_3}{W}=p^2/4$. This shows that 
	$$0 \le \I E\left((D-p/2)^2\right) = \I E(D^2) - p\,\I E(D) + p^2/4  = 0,
	$$
	  and therefore $D=p/2$ with probability~1.
\epf

Thus, by Lemma~\ref{lm:regular}, we have the following.

\begin{claim}\label{cor:addleafdensity}
	If $\h{H}{U}=\h{H}{W}$ for some graph $H$, then $\h{H'}{U}=\h{H'}{W}$ for any graph $H'$ that is obtained from $H$ by adding a pendant edge.
	In particular, $\h{F}{U}=\h{F}{W}$ for any forest~$F$.\qed
\end{claim}

We therefore obtain the following.

\begin{claim}\label{prop:3edgedensity} It holds that
	$\h{H}{U} = \h{H}{W}$ for each $H$ in $\C G_{5,3}=\{ K_{1,3},P_4,P_2\sqcup K_2,K_3\}$. In particular,  $\h{K_3}{U}=\h{K_3}{W}=0$ and thus $\h{H}{U}=0$ for every graph $H$ that contains a triangle. 
	
\end{claim}

\bpf 
	For the first three graphs of $\C G_{5,3}$, this follows immediately from Claim~\ref{cor:addleafdensity}.
	
	 Lemma~\ref{lm:Csoka}, when applied with $m=3$ to each of $U$ and $W$, gives the same linear relation (with all coefficients non-zero) relating the densities of the four graphs in~$\C G_{5,3}$. Since we have already established that $\h{H}{U}=\h{H}{W}$ for every $H\in \C G_{5,3}\setminus\{K_3\}$, we must also have that $\h{K_3}{U}=\h{K_3}{W}$.  
\epf 

Thus by Claims~\ref{cor:addleafdensity} and~\ref{prop:3edgedensity}, for each $m\in [10]$, any graph in $\C G_{5,m}$ that has different densities in $U$ and $W$ must belong to $\C H_m$, which we define to consist of $H\in \C G_{5,m}$ such that $H$ is triangle-free and $H$ is not a forest. Clearly, every such $H$ must have an induced cycle of length $4$ or $5$, which makes it easy to enumerate all graphs in~$\C H_m$.

We have that $\C H_4=\{C_4\}$ consists only of the 4-cycle and thus Lemma~\ref{lm:Csoka} gives that \beq{eq:C4}
\h{C_4}{U}=\h{C_4}{W}.
\eeq 

Next, $\C H_5$ consists only of the $5$-cycle $C_5$ and $C_4'$, a $4$-cycle with a leaf attached to one of its vertices. By~\eqref{eq:C4} and Claim~\ref{cor:addleafdensity}, we have that $\h{C_4'}{U}=\h{C_4'}{W}$. Thus, by Lemma~\ref{lm:Csoka}, the other graph in $\C H_5$, namely the $5$-cycle, must also have the same density in $U$ as in~$W$.

Since $\C H_6$ contains
only the complete bipartite graph $K_{2,3}$,  Lemma~\ref{lm:Csoka} gives that
$\h{K_{2,3}}{U}=\h{K_{2,3}}{W}$.

We now consider the random variable $\Pxy:=\codeg^U(x,y)=\int_0^1 U(x,z)U(z,y)\dd z$, the density of copies of $P_2$ which have $x,y$ as endpoints, where $x$ and $y$ are chosen uniformly and independently from~$[0,1]$. The following identities 
show that the first three moments of $\Pxy$ remain the same if we replace $U$ by~$W$:
\begin{eqnarray*}
	\I E(\Pxy)\ =\ \h{P_3}{U} &=& \h{P_3}{W}\ =\ p^2/4,\\
	\I E(\Pxy^2) \ =\ \h{C_4}{U}& =& \h{C_4}{W}\ =\ p^4/8,\\
	\I E(\Pxy^3) \ =\ \h{K_{2,3}}{U}& =& \h{K_{2,3}}{W}\ =\ p^6/16.
\end{eqnarray*}

We now observe that
$$
\I E\left(\Pxy(\Pxy-p^2/2)^2\right) = \I E(\Pxy^3) -p^2\, \I E(\Pxy^2) + \frac{p^4}4\, \I E(\Pxy) =\frac{p^6}{16}-p^2\,\frac{p^4}{8}+\frac{p^2}4\cdot \frac{p^2}4 = 0.
$$ 	
Since $\Pxy(\Pxy-p^2/2)^2 \ge 0$ deterministically, we have that $\Pxy\in \{0,p^2/2\}$ with probability~1. By $\I E(\Pxy)=p^2/4$, we conclude that
$$
	\pr(\Pxy=0)  = \pr(\Pxy=p^2/2) =\frac{1}{2}.
$$

Thus almost all pairs come in two types. Let $C$ consist of
those $(x,y)\in [0,1]^2$ for which $\codeg^U(x,y)=p^2/2$. Its complement consists a.e.\ of pairs with zero codegree. 
Since the measure of $C\subseteq [0,1]^2$ is $1/2$, we have $\int_0^1\deg_C(x)\dd x=1/2$, where $\deg_C(x)$ denotes the measure of 
$$
N_C(x):=\{y\mid (x,y)\in C\},\quad \mbox{for $x\in [0,1]$}.
$$

We know from Claim~\ref{cor:addleafdensity}
that $\h{P_5}{U}=\h{P_5}{W}$ where $P_5$ is the path with $5$ vertices.
One can compute $\h{P_5}{U}$ as follows. Recall that $P_5$ visits vertices $1,\dots,5$ in this order. First sample $x_1,x_3,x_5\in [0,1]$ and then pick common neighbours $x_2$ and $x_4$ of $x_1x_3$ and $x_3x_5$ respectively.
The measure of choices of $(x_2,x_4)\in [0,1]^2$ is $(p^2/4)^2$ if  $x_1x_3,x_3x_5\in C$ and 0 otherwise (apart from a null set of $(x_1,x_3,x_5)$).  The same argument applies to $\h{P_5}{W}$. Since $p\not=0$, the measure of $(x_1,x_3,x_5)\in [0,1]^3$ with $x_1x_3$ and $x_3x_5$ in $C$ must be the same as the analogous quantity for $W$, that is,~$1/4$. Thus, by 
the Cauchy-Schwarz Inequality, we have 
$$
 \frac14=\int_0^1 (\deg_C(x_3))^2\dd x_3\ge \left(\int_0^1 \deg_C(x_3)\dd x_3\right)^2=\frac14.
 $$
 We conclude that $\deg_C$ is the constant-$1/2$ function a.e.

For a.e.\ $x\in [0,1]$, the set $N_C(x)$ is independent in $U$. Indeed, 
the density of $C_5$ in $U$ can be written as the integral over $x_3\in [0,1]$ and then over $x_1,x_5\in N_C(x_3)$ of $U(x_1,x_5)(p^2/2)^2$. On the other hand, we know that $\h{C_5}{U}=\h{C_5}{W}=0$, giving the claim.

Pick a typical $x\in [0,1]$, that is, with $A':=N_C(x)$ being measurable, independent and of measure $1/2$. By above and the $(p/2)$-regularity of $U$, each vertex of $A'$ has degree $p/2$ in $U$ and almost all these edges connect $A'$ to its complement $B':=[0,1]\setminus A'$. Thus 
$$
 \frac p4=\int_{A'}\deg_{B'}^U(x)\dd x=\int_{A'\times B'} U(x,y)\dd x\dd y\le \frac12\,\h{K_2}{U},
 $$ 
where the factor $1/2$ arises because integrating over $B' \times A'$ would give exactly the same result, since $U$ is symmetric.
 But $\h{K_2}{U}=\h{K_2}{W}=p/2$, so these are all the edges, that is, $B'$ is an independent set.
 
By applying a measure-preserving transformation to $U$, we can assume that $A'=A$ and $B'=B$. These sets are independent in both $U$ and $W$. 
Since $W$ takes constant value $p$ on $A\times B$,
Lemma~\ref{lm:BipQR} gives that $U$ is also $p$ a.e. on $A\times B$, finishing the proof of
Theorem~\ref{th:1param}.

\section{Graphons with large density of independent $k$-sets}\label{containers}

We will need the following auxiliary result which, informally speaking, states that if a graph has many independent sets of large but fixed size $k$ then the graph has a large almost independent set. Let $\C I(G)$ denote the family of all independent sets in a graph $G$ and let $$
\C I_k(G):=\{I\in\C I(G)\mid |I|=k\}
$$ 
consist of all independent sets of size~$k$.

\begin{theorem}\label{th:ManyEk}
	For every $\delta>0$ there exists $\e>0$ such that for any $k\ge 1/\e$ there exists $n_0$ such that for every graph $G$ on $n\ge n_0$ vertices and every real $\alpha$, if $|\C I_k(G)|\ge (\alpha-\e)^k{n\choose k}$, then there exists $A\subseteq V(G)$ with $|A|\ge (\alpha-\delta)n$ and $e(G[A])\le \delta n^2$.
	\end{theorem}

\bpf Given $\delta>0$, choose sufficiently small $\e>0$; in particular, assume that $\e<\delta/2$. Given any $k\ge 1/\e$, let $n$ be sufficiently large. Let a graph $G=(V,E)$ and a real $\alpha$ satisfy the assumptions of the lemma. Assume that $\alpha>\delta$ as otherwise we can trivially let $A$ be the empty set.

We use a basic version of the container method that was introduced in high generality independently by 
Balogh, Morris and Samotij~\cite{BaloghMorrisSamotij15jams} and by Saxton and Thomason~\cite{SaxtonThomason15}, and whose roots go back to Kleitman and Winston~\cite{KleitmanWinston80,KleitmanWinston82} and Sapozhenko~\cite{Sapozhenko01,Sapozhenko02}.
Roughly speaking, we will encode each independent set $I$ of $G$ by a very small set $T\subseteq I$ together with a decoding procedure that produces a container $C=C(T)$ that necessarily contains $I$ and spans few edges. Then we take for $A$ a container $C(T(I))$ of the maximal size over all choices of~$I\in\C I_k(G)$. 

Formally, we proceed as follows. Assume that $V=[n]$ with the natural order.

Take any $I\in\C I(G)$. Enumerate $I=\{i_1<\dots<i_m\}$ using the natural order on $V=[n]$. The \emph{encoding} procedure produces $T=T(I)$ as follows.  Initially, let $T:=\emptyset$ and $j:=1$. Iterate the following step. 
Given the current values of $j\le m$ and $T\subseteq \{i_1,\dots,i_{j-1}\}$, add $i_j$ into $T$ if 
and only if
 \begin{equation}\label{eq:cond}
 |N(T\cup\{i_j\})|\ge |N(T)|+ \delta n/2,
 \end{equation}
 that is, $i_j$ has at least $\delta n/2$ neighbours outside of~$N(T)$.
 Then increase $j$ by $1$ and, if the new $j$ is still at most $m$, repeat the iteration step. 
 
Let $T=T(I)$ be the final set~$T$. Since only the vertices of $I$ were considered for inclusion into $T$, we have that $T\subseteq I$. Also, every time we add a vertex to $T$, the size of the neighbourhood $N(T)$ increases by at least $\delta n/2$. Thus
 \beq{eq:|T|}
  |T|\le 2/\delta.
  \eeq

Now, let us describe the \emph{decoding} procedure which constructs the \emph{container} $C(T)$ for any independent set $T\subseteq V$ of the graph~$G$. Let $t:=|T|$. Enumerate
$V\setminus T=\{v_1<\dots<v_{n-t}\}$, again using the natural order on~$[n]$. Initially, let $C:=T$ and $j:=1$. Repeat the following step given the current values of $j\le n-t$ and $C\subseteq T\cup \{v_1,\dots,v_{j-1}\}$: 
include $v_j$ into $C$ if and only if $v_j\not\in N(T)$ and 
 \begin{equation}\label{eq:C}
 |N(T_j\cup\{v_j\})|< |N(T_j)|+ \delta n/2, \quad \mbox{where }T_j:=\{v\in T\mid v<v_j\}.
 \end{equation}
 Then increase $j$ by 1 and, if $j\le n-t$, repeat the iteration step.
 
By construction, the final set $C=C(T)$ contains $T$ and is disjoint from $N(T)$. Also, in the notation of~\eqref{eq:C}, each vertex $v_j$ of $C\setminus T$ has fewer than $\delta n/2$ neighbours in $V\setminus N(T_j)$. Note that the last set contains $V\setminus N(T)\supseteq C$ (since $T_j\subseteq T$ and $C\cap N(T)=\emptyset$). Thus $C$ spans at most $|T|n+|C|\delta n/2$ edges. This is at most $2n/\delta+ \delta n^2/2<\delta n^2$ if $T$ satisfies~\eqref{eq:|T|}.

Let us also show that 
\beq{eq:ISubsetC}
 I\subseteq C(T(I)),\quad\mbox{for every independent set $I$ of $G$.}
\eeq
Let $T:=T(I)$. Since $T\subseteq C(T)$, it remains to  show that every $i_s\in I\setminus T$ belongs to $C(T)$, where $i_1<i_2<\ldots$ enumerate all elements of~$I$. The reason why $i_s$ was not included into $T$ when it was considered at the appropriate encoding step must be that~\eqref{eq:cond} fails for $j=s$, that is, $i_s$ adds fewer than $\delta n/2$ new neighbours when added to $\{ v\in T\mid v<i_s\}$. 
This is exactly the statement in~\eqref{eq:C}. Also, since $I\supseteq T$ is an independent set, we have that $i_s\not\in N(T)$. Thus $i_s\in C(T)$ by the definition of the decoding procedure, proving~\eqref{eq:ISubsetC} as desired.

By~\eqref{eq:|T|} and~\eqref{eq:ISubsetC}, we get the following upper bound on the number of independent sets of size exactly~$k$:
$$
|\C I_k(G)|\le \sum_{t=0}^{\floor{2/\delta}}
\sum_{T\in \C I_t(G)}
{|C(T)\setminus T|\choose k-t}.
$$ 
 Fix an index $t$, between $0$ and $\floor{2/\delta}$, such that the $t$-th summand is at least the average value, which is at least $1/(2/\delta+1)$ times~$|\C I_k(G)|$.  Given this $t$, let $A$ be a maximum-size container $C(T)$ over all independent $t$-sets $T\subseteq V$.
 Then
 \beq{eq:A}
\frac1{2/\delta+1} (\alpha-\e)^k{n\choose k}\le  {n\choose t} {|A|-t\choose k-t}.
\eeq

The set $A$, as some container $C(T)$ for a set $T$ with $|T|\le 2/\delta$, spans at most $\delta n^2$ edges in~$G$. Thus, in order to finish the proof of the theorem, we have to show that $|A|\ge (\alpha-\delta)n$. When $n$ tends to infinity (with $\delta\gg \e\ge 1/k$ fixed and $t\le 2/\delta$ bounded), the inequality in~\eqref{eq:A} gives that $|A|\to\infty$ and, in fact, $ |A|\ge  ((\delta/3)^{1/k}+o(1))c_tn$, where $c_t:=\left((\alpha-\e)^k/{k\choose t}\right)^{\frac{1}{k-t}}$. It is enough to show that e.g.\ $c_t\ge (\alpha-\delta/2)$. The last inequality is equivalent to
$$
 \left(\frac{\alpha-\e}{\alpha-\delta/2}\right)^k\ge {k\choose t} (\alpha-\delta/2)^{-t},
$$
 which holds for all large $k$ because the left-hand side grows exponentially in $k$, while the right-hand side grows at most polynomially. (Recall that $t\le 2/\delta$ is bounded.)

 Thus the set $A$ has all required properties.\epf

For $k\in\I N$, let 
$$\alpha_{k}(W):=\I P(X_k(W)=0)$$ be the probability that the $k$-sample $\I G(k,W)$ has no edges. We have 
$$
\alpha_{k+m}(W)\le \alpha_{k}(W)\, \alpha_{m}(W),\quad \mbox{for all $k,m\in\I N$},
$$ 
because the subgraphs of $\I G(k+m,W)$ spanned by the first $k$  and the last $m$ vertices are independent and distributed as $\I G(k,W)$ and $\I G(m,W)$ respectively.  Thus, by the Fekete Lemma, the limit
\begin{equation}\label{eq:lim}
	\alpha_\infty(W):= \lim_{k\to\infty} (\alpha_{k}(W))^{1/k}
\end{equation}
exists. Clearly, $\alpha_\infty(W)$ remains the same if we replace $W$ by any weakly isomorphic graphon.

In order to prove Theorem~\ref{th:alpha}, which states in the above notation that $\alpha(W)=\alpha_\infty(W)$,  we need to present some definitions and results related to measure theoretic aspects of graphons from~\Lo{Chapters 8 and 13}. Let $U$ and $W$ be graphons. We define the \emph{cut-norm} 
$$
\cutnorm{U}{W}:=\sup_{A,B\subseteq [0,1]} \left|\int_{A\times B}\left(U(x,y)-W(x,y)\right) \dd x\dd y\right|,
$$
where the supremum is taken over all pairs of measurable subsets of $[0,1]$.
For a measure-preserving function $\phi:[0,1]\to [0,1]$, 
the \emph{pull-back} $U^\phi$ of $U$ along $\phi$ is defined 
by 
$$
U^{\phi}(x,y):=U(\phi(x),\phi(y)),\quad x,y\in [0,1].
$$
It is routine to see that $U^{\phi}$ is a graphon which is weakly isomorphic to~$U$.
The \emph{cut-distance} is defined as
\beq{eq:CutDistance}
\delta_{\Box}(U,W):=\inf_{\phi} \cutnorm{U^\phi}{W},
\eeq
where the infimum is taken over all invertible measure-preserving maps $\phi:[0,1]\to [0,1]$.
See \Lo{Section~8.2} for more details and, in particular, \Lo{Theorem 8.13} for some alternative definitions that give the same distance. It can be easily verified that $\delta_\Box$ is a pseudo-metric on the space of graphons. Moreover, two graphons are weakly isomorphic if and only if they are at cut-distance~$0$, see e.g.\ \Lo{Theorem 13.10}.

\bpf[Proof of Theorem~\ref{th:alpha}] The inequality $\alpha_\infty(W)\ge \alpha(W)$ is easy. Indeed, 
pick an independent set $A\subseteq [0,1]$ in $W$ of  measure $\lambda(A)= \alpha(W)$ (which exists by~\cite[Lemma~2.4]{HladkyHuPiguet19})
and observe that the probability of seeing no edges in the $k$-sample $\I G(k,W)$ is at least $\lambda(A)^k$, the probability that all 
vertices land in~$A$.

Let us show the converse inequality $\alpha_\infty(W)\le \alpha(W)$. 
Let $\alpha:=\alpha_\infty(W)$ and assume that $\alpha>0$ as otherwise there is nothing to prove.
	Do the following for every $m\in\I N$. 
Let $\e>0$ be sufficiently small, in particular to satisfy Theorem~\ref{th:ManyEk} for $\delta:=1/m$.  By~\eqref{eq:lim}, pick $k\ge 1/\e$ such that $\alpha_{k}(W)\ge 4(\alpha-\e)^k$. Let $n$ be sufficiently large and
take the $n$-sample $G\sim\I G(n,W)$. Let $W_G$ be its graphon, that is, $W_G$ is the $n$-step $\{0,1\}$-valued graphon that encodes the adjacency relation of~$G$. 

It is easy to see that the mean of $\alpha_k(W_G)$ over $G\sim\I G(n,W)$ 
 is exactly $\alpha_k(W)$.
We claim that in fact $\alpha_k(W_G)$ is concentrated around this mean,
for which we apply Azuma's inequality (see e.g.~\cite[Theorem~2.25]{JansonLuczakRucinskiBook}).
Observe that 
$\alpha_k(W_G)$ can only change by at most $k/n$ if a vertex of $G$ is altered, i.e.\ 
the vertex-exposure martingale revealing $G$ and tracking $\alpha_k(W_G)$ is $(k/n)$-Lipschitz.
Setting $t:=2(\alpha-\e)^k$, Azuma's inequality states that
$$
\pr\left(\alpha_k(W_G) \le 2(\alpha-\e)^k\right) \le \pr \Big( \alpha_k(W_G) \le \I E(\alpha_k(W_G))-t\Big) \le \exp\left(\frac{-t^2}{2n (k/n)^2}\right) =o(1),
$$
where asymptotics are as $n\to \infty$.

It is also known that, as $n\to\infty$, the probability that the cut-distance between $W_G$ and $W$ is more than $o(1)$ is at most $o(1)$, specifically (see e.g.\ Lemma~10.16 in~\cite{Lovasz:lngl})
$$
\I P\left(\,\delta_\Box(W_G,W)>22/\sqrt{\log n}\,\right)\le  \exp(-n/(2\log n)).
$$
Thus, for large enough $n$, there is a graph $G$ on $[n]$ with $\alpha_k(W_G)\ge 2(\alpha-\e)^k$ and $\delta_\Box(W,W_G)\le \delta/2$ because $\I G(n,W)$ satisfies each of these properties with probability $1-o(1)$ as $n\to\infty$. Since $|\C I_k(G)|/{n\choose k}=\alpha_k(W_G)+o(1)$, Theorem~\ref{th:ManyEk} applies to the graph $G$ for large enough $n$ and returns a set $A'$ of vertices of size at least $(\alpha-\delta)n$ spanning at most $\delta n^2$ edges. Let $A:=\cup_{i\in A'} [\frac{i-1}n,\frac in)$ be the subset of $[0,1]$ corresponding to $A'$ when we pass from $G$ to its graphon~$W_G$. Take an invertible measure-preserving map $\phi:[0,1]\to[0,1]$ such that $\cutnorm{W^\phi}{W_G}<\delta$. Defining $S_m:=\phi(A)$ to be the image of the set $A$, we obtain that 
$$
\int_{S_m^2} W(x,y)\dd x\dd y=\int_{A^2} W^\phi (x,y)\dd x\dd y\le \int_{A^2} W_G(x,y)\dd x\dd y+\cutnorm{W^\phi}{W_G}\le 2\delta= 2/m.
$$ 
 
We now proceed as in \cite{HladkyHuPiguet19}. Recall that a sequence of functions $f_1,f_2,\dots$ in $L^\infty([0,1],\lambda)$, the dual space of $L^1([0,1],\lambda)$, \emph{weak-$*$ converges} to $f\in L^\infty([0,1],\lambda)$ if 
\beq{eq:Weak*}
 \lim_{n\to\infty} \int_0^1 f_n(x)g(x)\dd x=\int_0^1 f(x)g(x)\dd x,\quad\mbox{for  every $g\in L^1([0,1],\lambda)$.}
 \eeq
By the Sequential Banach-Alaoglu Theorem 
(see e.g.\ \cite[Theorem~1.9.14]{Tao10erira}),
the sequence of the characteristic functions of the sets $S_m$ viewed as elements of $L^\infty([0,1],\lambda)$ has a subsequence that weak-$*$ converges to some function~$f\in L^\infty([0,1],\lambda)$. Note that $f(x)\ge 0$ for $\lambda$-a.e.\ $x\in [0,1]$. Indeed, letting $g:=\chi_{X}$ be the characteristic function of the measurable set $X:=\{x\in [0,1]\mid f(x)<0\}$, we get from~\eqref{eq:Weak*} that
 $$
 0\ge \int_X f\dd \lambda=\int_0^1 fg\dd \lambda=\lim_{m\to\infty} \int_0^1 \chi_{S_m}\,g\dd\lambda\ge 0,
$$
 from which it follows that $\lambda(X)=0$. Likewise, we obtain that $f\le1$ a.e.\ on~$[0,1]$.

Let the \emph{support} of a function $g:[0,1]\to \I R$ be the set $\supp(g):=\{x\in [0,1]\mid g(x)\not=0\}$.
Lemma~2.4 in \cite{HladkyHuPiguet19} states in fact that, for any graphon $W$, the support of the weak-$*$ limit of the characteristic functions of $W$-independent sets is $W$-independent. Thus  $S:=\supp(f)$, the support of $f$, is an independent set in~$W$. By the definition of weak-$*$ convergence and the fact that $\|f\|_\infty\le
1$, we have that 
$$
 \alpha(W)\ge \lambda(S)\ge \int_0^1 f(x)\dd x=\lim_{m\to\infty}\int_0^1 \chi_{S_m}(x)\dd x=\lim_{m\to\infty} \lambda(S_m)=\alpha_\infty(W).
 $$
  
 This shows that $\alpha(W)=\alpha_\infty(W)$, proving Theorem~\ref{th:alpha}.\epf

\section{Proof of Theorem~\ref{th:Negated}}\label{sec:Negated}
\renewcommand{\beta}{p}
\renewcommand{\gamma}{q}

\bpf[Proof of Theorem~\ref{th:Negated}] Recall that $W$ is the 2-step graphon with steps $A$ and $B$ which is $0$ on $A^2$, $\beta$ on $A\times B$ and $\gamma$ on~$B^2$, and where $A=[0,a)$. 
Let $U$ be an arbitrary graphon such that for every $k\in \I N$ the distributions of $X_k(U)$ and $X_k(W)$ are the same; let us denote this random variable by~$X_k$. We have to show that $U$ is weakly isomorphic to~$W$. Assume that $a\in(0,1)$ as otherwise $W$ is weakly isomorphic to a 1-step graphon and the conclusion follows from the results in~\cite{Csoka16}.
 
By Theorem~\ref{th:alpha} we conclude that $\alpha(U)=\lim_{k\to\infty} (\I P(X_k=0))^{1/k}=\alpha(W)$. Thus, by~\cite[Lemma~2.4]{HladkyHuPiguet19}, there is a set $A'$ of measure $a$ with
$U$ being 0 on $A'\times A'$ a.e. By taking a measure-preserving Borel isomorphism $\phi:[0,1]\to [0,1]$ with $\phi(A')=A$ and replacing  $U$ with $U^\phi$ (which is weakly isomorphic to $U$), we can assume that $A'=A$. Recall that $\deg_A^U(x):=\int_A U(x,y)\dd y$ for $x\in [0,1]$.

\begin{claim}\label{cl:dA} For almost every $x\in B$, we have $\deg_A^U(x)\ge a \beta$.
\end{claim}
 \bcpf If the claim is false, then by the continuity of measure, there is $\e>0$ such that the measure of 
 $B':=\{x\in B\mid \deg_A^U(x)\le a\beta-\e\}$ is at least $\e$. Take $k$ sufficiently large. Let us lower bound $\alpha_{k}(U)$, the probability that
 the $k$-sample $\I G(k,U)$ spans no edges. 
 Recall that we sample uniform $x_1,\dots,x_k\in [0,1]$ and then make each pair $ij$ an edge with probability $U(x_i,x_j)$, with all choices being mutually independent. With probability $a^k$, all elements $x_i$ belong to $A$ (when almost surely $\I G(k,U)$ has no edges). A disjoint event is that $x_1$ belongs to $B'\subseteq B$, which has probability is at least $\e>0$.
 Conditioned on this event, the probability of having no edges is
 at least 
 $$
  \left(\int_A (1-U(x_1,y))\dd y\right)^{k-1}=(a-\deg^U_A(x_1))^{k-1}\ge (a(1- \beta)+\e)^{k-1},
$$
 since, when ignoring null sets, it is enough that all other $k-1$ vertices belong to $A$ and are all non-adjacent to~$x_1$.

Thus $\alpha_{k}(U)\ge a^k+\e (a(1-\beta)+\e)^{k-1}$. For $W$, it is easy to write an explicit formula, where~$i$ denotes the number
 of sampled vertices that belong to~$B$:
 \begin{eqnarray*}
 \alpha_{k}(W)&=&\sum_{i=0}^k {k\choose i} (1-a)^i a^{k-i} (1-\beta)^{(k-i)i} (1-\gamma)^{{i\choose 2}}\\
 &=&a^k + k(1-a)a^{k-1} (1-\beta)^{k-1}+\dots\ .
 \end{eqnarray*}
Note that the first term $a^k$ matches that for~$U$. 
Of course, we have $\alpha_k(U)=\I P(X_k=0)=\alpha_k(W)$. 
Thus a desired contradiction, namely that $\alpha_{k}(U)>\alpha_{k}(W)$, will follow if we show that for every $i\in [k]$,
 \beq{eq:aim1}
 \e (a(1-\beta)+\e)^{k-1}> k\cdot {k\choose i} (1-a)^i  a^{k-i}(1-\beta)^{(k-i)i} (1-\gamma)^{{i\choose 2}}.
 \eeq
 Informally speaking, if $i$ is small, then the main terms are $(a(1-\beta)+\e)^k$ versus $(a(1-\beta)^i)^k$; otherwise either $(1-\beta)^{(k-i)i}$ or $(1-\gamma)^{{i\choose 2}}$ (and thus the right-hand side of~\eqref{eq:aim1}) is very small.
 Formally, given $\e>0$ as above, fix a large constant $M\gg 1/\e$ and then let $k\to\infty$. If $1\le i\le M$  then  the ratio of the right-hand side to the left-hand side of~\eqref{eq:aim1} is at most 
 $$
 O\left(\left(\frac{a(1- \beta)}{a(1-\beta)+\e}\right)^k\cdot k^{M+1}\right)=o(1).
 $$ 
 If $M< i\le k$, then since (slightly crudely) $\max\left((k-i)i,{i\choose 2}\right)\ge ki/4$ the ratio is at most
 $$
 \max(1-\beta,1-\gamma)^{ik/4}\cdot \frac{k^{i+1}}{\e (a(1-\beta)+\e)^k}=o(1),
 $$
where the last estimate holds since $p,q \neq 0$ by assumption.
  This proves~\eqref{eq:aim1} for all large $k$ and finishes the proof of the claim. 
 \ecpf

Let $U'$ be the graphon obtained from $U$ by averaging it over $(A\times B)\cup (B\times A)$ and over~$B^2$. That is, $U'$ is the 2-step graphon with parts $A$ and $B$ which assumes value 0 on $A^2$, value $\gamma':=\frac1{(1-a)^2}\,\int_{B^2} U(x,y)\dd x\dd y$ on $B^2$, and value 
\beq{eq:beta'}
 \beta'
 :=\frac1{a(1-a)} \int_{A\times B} U(x,y)\dd x\dd y = \frac1{a(1-a)} \int_{B} \deg^U_A(y)\dd y\ge \beta
 \eeq
 on $(A\times B)\cup (B\times A)$, where we applied 
 Claim~\ref{cl:dA} in~\eqref{eq:beta'}.
Consider the density of $P_3$, the path visiting vertices $1,2,3$ in this order. Its density, say in $U$, can be written as $\h{P_3}{U}=\int (\deg^U(x_2))^2\dd x_2$. 
Clearly, when we pass from $U$ to $U'$ then the degrees inside $A$ (resp.\ degrees inside $B$) are all replaced by their average value. Thus the average of $(\deg(x))^2$ over $x$ in $A$ (resp.\ $B$) does not increase by the Cauchy--Schwarz Inequality. By adding up these two averages weighted by $a$ and $1-a$ respectively, we get the average of $(\deg(x))^2$ over $x\in [0,1]$, which is the density of~$P_3$. Thus $\h{P_3}{U}\ge \h{P_3}{U'}$. 
The $W$-degrees of $x$ in $A$ and $B$ are the constants $(1-a)\beta$ and $a\beta+(1-a)\gamma$ respectively. These constants for $U'$ are $(1-a)\beta'$ and $a\beta'+(1-a)\gamma'$ respectively. 
Since $\beta'\ge\beta$ by~\eqref{eq:beta'}, we have that $(1-a)\beta'\ge (1-a)\beta$.
Furthermore, since $U'$ and $W$ have the same edge density (which can be computed as the convex combination of the average degrees in $A$ and $B$ weighted by $a$ and $1-a$), we have that $a\beta'+(1-a)\gamma'\le  a\beta+(1-a)\gamma$. 
Thus,
$$
 (1-a)\beta'\ge (1-a)\beta\ge a\beta +(1-a)\gamma\ge a\beta'+(1-a)\gamma',
$$
 where the middle inequality is an assumption of the theorem.
 We see that when we
pass from $U'$ to $W$, the degrees get more even (with the average  staying the same)
and so the density of $P_3$, the average of $(\deg(x))^2$, does not increase. 
However, the density of $P_3$ is determined by $X_3$ by Lemma~\ref{lm:2Edge}. Thus the degree functions of $U$, $U'$ and $W$ coincide a.e.  Thus $\beta'=\beta$, $\gamma'=\gamma$, 
\beq{eq:DegInB}
 \deg_A^U(x)=a\beta,\
 \deg_B^U(x)=(1-a)\gamma \mbox{\ for a.e.\ $x\in B$\quad and\quad $\deg_B^U(y)=(1-a)\beta$\ for  a.e.\ $y\in A$}.
 \eeq

\begin{claim}\label{cl:codegree} For almost every $(x,y)\in B^2$, we have that $\codeg_A^U(x,y)
	= a\beta^2$. (Recall that we denote $\codeg_A^U(x,y):=\int_A U(x,z)U(y,z)\dd z$.)\end{claim}

\bcpf Suppose first that $\codeg_A^U(x,y)> a\beta^2$ for some set of $(x,y)\in B^2$ of positive measure.
Then by the continuity of measure there exists $\e>0$ such that
the measure of 
$$
B':=\left\{(x,y)\in B^2\mid \codeg_A^U(x,y)\ge a\beta^2+\e\right\}
$$
 is at least~$\e$.

When we compute $\I P(X_k\le 1)$ using $U$ or $W$, the $k$-tuples of vertices that have at most one point in $B$ contribute the same amount. For example, let us condition on $x_i$ being the unique sampled vertex that belongs to $B$. Then  each other $x_\ell\in A$ is adjacent ot $x_i$ with probability $\frac1a\deg_A^U(x_i)$ which is equal to $\frac1a\deg_A^W(x_i)$ by~\eqref{eq:DegInB}. Moreover, these choices for different choices of $\ell$ are independent of each other, both in $U$ and in $W$. Thus any particular adjacency pattern of $x_i\in B$ to the other $k-1$ vertices from $A$ has the same conditional probabilities in $U$ and in~$W$. 

Consider the remaining contribution to $\I P(X_k\le 1)$, i.e.\ when at least two sampled vertices belong to~$B$. First, take~$U$. With probability at least $\e$, the pair $(x_1,x_2)$ belongs to the set $B'$. Conditioned on this, each other vertex $x_\ell$ is adjacent to neither $x_1$ nor $x_2$ with probability
 $$ 
 \int_A (1-U(x_1,y))(1-U(x_2,y))\dd y = a-2a\beta +\codeg_A^U(x_1,x_2),
 $$ 
with these choices being mutually independent for different values of~$\ell$. This contributes at least $\e (a(1-\beta)^2+\e)^{k-2}$ to~$\I P(X_k\le 1)$. On the other hand, an explicit summation formula can be written for $W$. First assuming that $\beta,\gamma \neq 1$, we have by above that
\begin{align*}
\e (a(1-\beta)^2+\e)^{k-2} & \le \sum_{i=2}^{k} {k\choose i} a^{k-i} (1-a)^i 
(1-\beta)^{(k-i)i}(1-\gamma)^{{i\choose 2}}\left(1+
{i\choose 2}\frac{\gamma}{1-\gamma}
+(k-i)i\frac{\beta}{1-\beta}\right).
\end{align*}
As before by taking $k\to\infty$ and looking at the cases  $i=O(1)$ and $i\gg 1$ separately, one can argue that $\e (a(1-\beta)^2+\e)^{k-2}$ is strictly larger than $k$ times the maximum term in the sum, giving a
contradiction.
If either $\beta$ or $\gamma$ is $1$, the inequality must be rewritten to avoid dividing by zero, and in fact the only
non-zero terms come from $i=2$ (if $\gamma=1$) or $i=k$ (if $\beta=1$).
It is then easy to obtain the necessary contradiction directly.

Thus $\codeg_A^U(x,y)\le a\beta^2$ for a.e.~$(x,y)\in B^2$.
The integral of $\codeg_{A}(x,y)$ over all $(x,y)\in B^2$ can be written 
as the integral over $z\in A$ of $(\deg_B(z))^2$. By~\eqref{eq:DegInB}, the latter integral is the same for $U$ as for~$W$. Thus 
$$
 (1-a)^2 a\beta^2 \ge \int_{B^2}\codeg_A^U(x,y)\dd x\dd y=\int_{B^2}\codeg_A^W(x,y)\dd x\dd y=(1-a)^2 a\beta^2,
$$
 and the first integrand must be $ap^2$ a.e.\ on~$B^2$, proving the claim.\ecpf

It follows that the density of triangles with one vertex in $A$ and two vertices in $B$ is the same in $U$ as in $W$. Indeed, first sample two vertices $x,y$ in $B$, which are adjacent in $U$ and in $W$ with the same probability $\gamma$ by~\eqref{eq:DegInB}, and observe that the probability of a vertex from $A$ being adjacent to both is exactly $\frac1a\,\codeg_A(x,y)$. Also, $U$ and $W$ have zero density of triangles  with at least two vertices in $A$.
Since $U$ and $W$ have the same triangle density, namely $\I P(X_3=3)$, they must have the same density of triangles that lie inside~$B$. 

Thus the density of triangles with any given partition of their vertices between $A$ and $B$ is the same for $U$ as for~$W$. Since the degrees in $U$ are the same as the degrees in $W$ by~\eqref{eq:DegInB}, this allows us to conclude by a version of Lemma~\ref{lm:regular} that the density of $K_3'$, the triangle with a pendant edge, is the same in $U$ as in~$W$. (Indeed, this is true even if we specify, relative to $A$ and $B$, where the vertices of the triangle in $K_3'$ lie.)  Lemma~\ref{lm:Csoka} applied to $\C G_{4,4}=\{C_4,K_3'\}$ gives that $U$ and $W$ have the same density of 4-cycles.

Claim~\ref{cl:codegree} implies that  the density of the 4-cycle (which we assume to visit the vertices $1,2,3,4$ in this order) conditioned on $x_1,x_3\in A$ and $x_2,x_4\in B$ is $\beta^4$, the fourth power of edge density~$\beta$ between $A$ and $B$. (Indeed, to sample such a $4$-cycle, we can first sample uniform $(x_2,x_4)\in B^2$; then we independently sample two vertices of $A$, each being adjacent to both $x_2$ and $x_4$ with probability $(\frac1a\, \codeg_A(x_2,x_4))^2$.) Thus if we let $U'$ be obtained from $U$ by letting it be 0 on $B^2$, then all assumptions of Lemma~\ref{lm:BipQR} are satisfied. The lemma gives that $U'$ (and thus $U$) assumes the constant value $\beta$ between $A$ and~$B$.

Thus we know all about $U$ except its values on~$B^2$. Our knowledge about $U$ is enough to compute the density of all types of $4$-cycles except those entirely inside~$B$. For example, if the sampled 4-cycle is to visit parts $B,B,B,A$ in this order, then we can first sample a 3-vertex path $x_1,x_2,x_3$ in $B$ (knowing its density since $\deg^U_B(x)=\deg^W_B(x)=(1-a)\gamma$ for almost every $x\in B$) and then use Claim~\ref{cl:codegree} to see that, conditioned on $x_4\in A$, the probability of $x_4$ being adjacent to both $x_1$ and $x_3$ is exactly $\beta^2$ (same as in~$W$). Since 
the graphons $U$ and $W$ have the same overall density of 4-cycles, they have the same density of 4-cycles inside~$B$. This relative density of 4-cycles is the fourth power of the relative edge density since $W$ is constant on~$B$.
By the Chung--Graham--Wilson Theorem (Theorem~\ref{th:CGW}),
the graphon $U$ assumes the constant value $\gamma$ on~$B^2$. We conclude that $U=W$ a.e., proving Theorem~\ref{th:Negated}.
 \epf

\section{Proof of Theorem~\ref{th:01p}}\label{01p}

\bpf[Proof of Theorem~\ref{th:01p}] Recall that $W$ is the 2-step graphon  with the first step $A=[0,a)$ being an independent set, while $W$ assumes value 1 on $B^2$ for $B:=[a,1]$ and value~$p$ on~$A\times B$. We have to show that if a graphon $U$ satisfies $\alpha(U)= a$, $\omega(U)= 1-a$ and $X_4(U)=X_4(W)=:X_4$, then $U$ is weakly isomorphic to~$W$.

By Theorem~\ref{th:alpha} (and~\cite[Lemma~2.4]{HladkyHuPiguet19}), there are subsets $C,D\subseteq [0,1]$ of measures $a$ and $1-a$ respectively such that $U$ is 0 on $C^2$ a.e.\ and $U$ is $1$ on $D^2$ a.e. The intersection $C\cap D$ must have measure 0, so we can assume that $C$ and $D$ partition~$[0,1]$. By applying a measure-preserving transformation to $U$, assume that $C=A$ and~$D=B$. 

\begin{claim}\label{cl:3} Almost every $x\in A$ satisfies  $\deg^U_B(x)=(1-a)p$.
\end{claim}
\bcpf The graphon $U$ has the same $K_4$-density as $W$ by $\h{K_4}{U}=\I P(X_4=6)$. 
Since $A$ is an independent set, there are only two types of $K_4$ in $U$: those inside $B$ (and their contribution to the overall density is $(1-a)^4$, the same as the analogous quantity for $W$) and those that have three vertices in $B$ and one vertex in~$A$. 
Since $U=1$ on $B^2$, the latter type of $4$-cliques determines $\frac1a\int_A (\deg_B^U(x))^3\dd x$, the third moment of the random variable $Y:=\deg_B^U(x)$, where $x$ is a uniform element of $A$.  This third moment is the same as for $W$, which is $((1-a)p)^3$ as $\deg_B^W(x)$ is the constant function $(1-a)p$. Also, we have by Lemma~\ref{lm:2Edge}  that
\beq{eq:IY}
 \I E(Y)=\frac1a\int_{A}\deg_B(x)\dd x=
\frac{\frac12\,\I P(X_2=1)-(1-a)^2}{a}=(1-a)p.
\eeq
 Thus $(\I E(Y))^3=\I E(Y^3)$ which for a non-negative variable $Y$ is possible only if $Y$ is constant a.e. Of course, the constant value of $Y$ must be $\I E(Y)=(1-a)p$, proving the claim.\ecpf 

Consider the random variable $Z:=\deg_A^U(x)$ for uniform $x\in B$. A calculation analogous to that in~\eqref{eq:IY} shows that $\I E(Z)=ap$. Our knowledge about $U$ directly gives the density of all possible copies of $P_3$, except $ABA$-paths (that is, copies of $P_3$ that have the middle vertex in $B$ and the other two in $A$). For example, the density of $BAB$-paths can be computed by sampling $x_2\in A$ first and then using Claim~\ref{cl:3}. Since the total $P_3$-density
in $U$ is the same as that for $W$
by Lemma~\ref{lm:2Edge}, we conclude that the density of $ABA$-paths in $U$ is also the same as in $W$, which is $(1-b)(ap)^2$. Thus $\I E(Z^2)=(ap)^2=(\I E(Z))^2$. This implies that for a.e.\ $x\in B$ we have $\deg_A^U(x)=ap$. (Alternatively, this conclusion can be reached by applying the argument of Claim~\ref{cl:3} to the complementary graphon~$1-U$.)

Let $K_4^-$ be the 4-clique minus an edge, the unique graph on $4$ vertices with $5$ edges.
Clearly, we have $\h{K_4^-}{U}=\frac16\, \I P(X_4=5) + \I P(X_4=6)=\h{K_4^-}{W}$.
The graphon $U$ has 2 types of $K_4^-$ of positive density. The first type consists of those copies of $K_4^-$ that have exactly 3 vertices in $B$ (and since $U$ is 1 a.e.\ on $B^2$, the corresponding density is determined by the degree distribution on $A$, which we know by Claim~\ref{cl:3}). Thus we know the density in $U$ of the other copies of $K_4^-$ which have 2 points in $A$ and 2 points in $B$ (and this matches that for $W$). Thus $U$ has the same density of $ABAB$-cycles as $W$. Also, their densities of $AB$-edges coincide, e.g.\ by Claim~\ref{cl:3}. Since $W$ is constant-$p$ on $A\times B$, the same must hold for $U$ by Lemma~\ref{lm:BipQR}. Thus $U$ and $W$ are weakly isomorphic.\epf

\section{Proofs of Propositions~\ref{pr:cycles} and~\ref{pr:Kkl}}\label{OtherQns}

\bpf[Proof of Proposition~\ref{pr:cycles}] We have to show that the family of graphs with at most one cycle is not forcing.
Here we use the observation that if $W$ is an $n$-step graphon with parts of measure $1/n$ and its values are encoded by a  symmetric $n\times n$ matrix $A\in [0,1]^{n\times n}$, then 
\beq{eq:Ck}
\h{C_k}{W}
=\frac1{n^k}\,\sum_{i=1}^n \lambda_i^k,\quad \mbox{for every $k\ge 3$,}
\eeq
where $\lambda_1,\dots,\lambda_n$ are the eigenvalues of $A$, repeated with their multiplicities. 
Indeed, $\h{C_k}{W}$ is  the sum over all ordered $k$-tuples $(v_0,\dots,v_{k-1})\in [n]^k$ of $\frac1{n^k}\,\prod_{i=0}^{k-1} A_{v_i,v_{i+1}}$, (where $v_k:=v_0$). On the other hand, the $j$-th diagonal entry of $A^k$ is the sum  of $\prod_{i=0}^{k-1} A_{v_i,v_{i+1}}$ taken over all $(v_0,\dots,v_{k-1})\in [n]^k$ with $v_0=j$. Summing this over all $j\in [n]$, we get that $\h{C_k}{W}$ is the trace of $\frac1{n^k}\,A^k$, giving the identity in~\eqref{eq:Ck}.

Take, for example,  the following unit vectors 
$$
\V x_1:=\frac1{\sqrt{3}}\Matrix{c}{1\\ 1\\ 1},\quad 
\V x_2:=\frac{1}{\sqrt2}\Matrix{r}{1\\ -1\\ 0},\quad\mbox{and}\quad
\V x_3:=\frac1{\sqrt6}\Matrix{r}{2\\ -1\\ -1},
$$
and let e.g.\ $\e:=1/4$. Note that $\V x_2$ and $\V x_3$ are orthogonal to~$\V x_1$. It routinely follows that the  symmetric $3\times 3$ matrices
\begin{equation}\label{eq:AA'}
	A:=\V x_1\V x_1^T+\e\, \V x_2\V x_2^T\quad\mbox{and}\quad A':=\V x_1\V x_1^T+\e\,\V x_3\V x_3^T
\end{equation} 
have the same eigenvalues (namely $1$, $\e$ and $0$), all entries in~$[0,1]$ and all row sums equal (namely, to $1$, which is the row sum of $\V x_1\V x_1^T$).

Let $W$ and $W'$ be the $3$-step graphons, with steps of measure $1/3$, whose values are given by the symmetric matrices $A,A'\in [0,1]^{3\times 3}$. The maximum entry of $A'$ is $1/3+2\e/3$, which is strictly larger than the maximum entry $1/3+\e/2$ of $A$, so the graphons $W$ and $W'$ are are not weakly isomorphic by e.g.\ considering the density of $K_k$ as $k\to\infty$ (and noting that the maximum entry of $A'$ is on the diagonal).

On the other hand, $W$ and $W'$ have the same cycle densities by~\eqref{eq:Ck}. Since they are both $(1/{3})$-regular, they have the same homomorphism density for every graph with at most one cycle by Lemma~\ref{lm:regular}.
This proves Proposition~\ref{pr:cycles}.\epf

\bpf[Proofs of Propositions~\ref{pr:Kkl}] We have to show that the family of graphs of diameter at most $d$ is not forcing. Let
$W$ (resp.\ $W'$) be the graphon of the disjoint union $G:=P_{d+2}\sqcup P_{d+2}$ (resp.\ $G':=P_{d+3}\sqcup P_{d+1}$). (Recall that $P_n$ denotes the path with $n$ vertices.) In other words, each of $W$ and $W'$ is the  $\{0,1\}$-valued  step graphon with $2d+4$ steps of equal measure that encodes the adjacency relation of the corresponding graph. They are not weakly isomorphic
because the induced density of $P_{d+3}$ is zero in $W$ but not in~$W'$.

On the other hand, $\h{F}{W}=\h{F}{W'}$ for every graph $F$ of diameter at most~$d$. Indeed, if $V(F)=[k]$, then for each choice of $\V x=(x_1,\dots,x_k)\in [0,1]^k$ for which the integrand in~\eqref{eq:t=Int} is positive, the union of the corresponding vertices of $G$ or $G'$ must induce a subgraph of diameter at most $d$, that is, a sub-path with $i\le d+1$ vertices. The number of sub-paths of any given order $i\le d+1$ is the same for $G$ and $G'$ (namely, $2d+4-2i$). Thus $\h{F}{W}=\h{F}{W'}$.

We conclude that the family of graphs of diameter at most~$d$ is not forcing.\epf

\begin{remark} One can view the construction in the proof of Proposition~\ref{pr:Kkl} as first taking the graphon $W_{C_m}$ of the $m$-cycle $C_m$ for $m:=2d+4$ and then decreasing density to zero on some two edges. Instead, we could have first multiplied the whole $m$-cycle graphon $W_{C_m}$ by some real $p\in (0,1)$ and then defined $W$ and $W'$ by modifying $pW_{C_m}$ on the the same pairs of edges of $C_m$ as before but in a way such that the new graphons are still $(2p/m)$-regular but not weakly isomorphic. This modified construction shows  by Lemma~\ref{lm:regular} that we can increase the family in Proposition~\ref{pr:Kkl} by taking all graphs of diameter at most $d$ and then attaching any number of pendant trees.	
\end{remark}

\section{Concluding remarks}\label{sec:concluding}

The only information that our results on Question~\ref{q:Sos} used was the conclusions about the densities for graphs on at most $5$ vertices that follow from the distribution of $X_2,\ldots,X_5$, and the probabilities that $X_k$ is $0$, $1$, ${k\choose 2}-1$, or ${k\choose 2}$
as $k\to\infty$. It would be interesting to find new types of arguments that use much more substantial information about the random variables~$X_k$.

An intriguing open question (which is a considerable weakening of Question~\ref{q:Sos}) is whether
the sequence $(X_k(W))_{k\in\I N}$ determines the essential supremum $\|W\|_\infty$
of an arbitrary graphon~$W$. If true, this would greatly enlarge the set of 2-step graphons $W$ for which we can prove that $W \in \SosFamily$.
 One can show that the limit of $\left(\I P\left(X_k(W)={k\choose 2}\right)\right)^{-{k\choose 2}}$ as $k\to\infty$ is the supremum of $\int_{A^2} W(x,y)\dd x\dd y$ over all subsets $A\subseteq [0,1]$ of positive measure. However this ``symmetric'' supremum need not be equal to $\|W\|_\infty$ (take e.g.\ a 2-step bipartite graphon of Theorem~\ref{th:1param}).

\bibliographystyle{plain}
\bibliography{SosQuestionAMH}

\end{document}